\documentclass[12pt,reqno]{article}
\usepackage{amsmath,amssymb,amsfonts,amscd,latexsym,amsthm,mathrsfs}
\theoremstyle{remark}

\let\wh\widehat

\renewcommand{\d}{{\mathrm d}}

\renewcommand{\Re}{\operatorname{Re}}
\begin{document}

\title{Transformations of $L$-values}

\author{Wadim Zudilin\thanks{This work is supported by Australian Research Council grant DP110104419.
The text is loosely based on my talk ``Mahler measures and $L$-series of elliptic curves''
at the conference ``Analytic number
theory\,---\,related multiple aspects of arithmetic functions''
(Research Institute for Mathematical Sciences, Kyoto University,
Japan, October~31--November~2, 2011).}\\[1mm]{\small School of Mathematical and Physical Sciences,}\\[-1.2mm]
{\small The University of Newcastle, Callaghan, NSW 2308, Australia}}

\date{February 2012}

\maketitle

\begin{abstract}
In our recent work with M.~Rogers on resolving some Boyd's conjectures on two-variate Mahler measures,
a new analytical machinery was introduced to write the values $L(E,2)$ of $L$-series
of elliptic curves as periods in the sense of Kontsevich and Zagier. Here we outline,
in slightly more general settings, the novelty of our method with Rogers,
and provide a simple illustrative example.
\end{abstract}


Throughout the note we keep the notation $q=e^{2\pi i\tau}$ for $\tau$
from the upper half-plane $\Re\tau>0$, so that $|q|<1$. Our basic constructor
of modular forms and functions is Dedekind's eta-function
\begin{equation*}
\eta(\tau):=q^{1/24}\prod_{m=1}^{\infty}(1-q^m)
=\sum_{n=-\infty}^{\infty}(-1)^nq^{(6n+1)^2/24}
\end{equation*}
with is modular involution
\begin{equation}
\eta(-1/\tau)=\sqrt{-i\tau}\eta(\tau).
\label{k01}
\end{equation}
We also set $\eta_k:=\eta(k\tau)$ for short.

We first describe a part of the general machinery from our joint works \cite{RZ10,RZ11} with M.~Rogers
on an example of computing the value $L(E_{32},2)$ of the $L$-series associated
with a conductor 32 elliptic curve. It is known~\cite{MO97} that the corresponding cusp form
in this case is $f_{32}(\tau):=\eta_4^2\eta_8^2$, so that
$L(E_{32},s)=L(f_{32},s)$. We choose the conductor 32 case here because it is not discussed in \cite{RZ10,RZ11}.

Note the (Lambert series) expansion
\begin{equation}
\frac{\eta_8^4}{\eta_4^2}
=\sum_{m\ge1}\biggl(\frac{-4}m\biggr)\frac{q^m}{1-q^{2m}}
=\sum_{\substack{m,n\ge1\\\text{$n$ odd}}}\biggl(\frac{-4}m\biggr)q^{mn},
\label{k02}
\end{equation}
where $\bigl(\frac{-4}m\bigr)$ is the quadratic residue character modulo~4.
In notation $\delta_{2\mid n}=1$ if $2\mid n$ and 0 if $n$~is odd,
we can write \eqref{k02} as
\begin{equation*}
\frac{\eta_8^4}{\eta_4^2}
=\sum_{m,n\ge1}a(m)b(n)q^{mn},
\qquad\text{where}\quad
a(m):=\biggl(\frac{-4}m\biggr), \quad b(n):=1-\delta_{2\mid n}.
\end{equation*}

Then
\begin{align*}
f_{32}(it)
&=\frac{\eta_8^4}{\eta_4^2}\,\frac{\eta_4^4}{\eta_8^2}\bigg|_{\tau=it}
=\frac{\eta_8^4}{\eta_4^2}\bigg|_{\tau=it}
\cdot\frac1{2t}\,\frac{\eta_8^4}{\eta_4^2}\bigg|_{\tau=i/(32t)}
\\
&=\frac1{2t}\sum_{m_1,n_1\ge1}a(m_1)b(n_1)e^{-2\pi m_1n_1t}
\sum_{m_2,n_2\ge1}b(m_2)a(n_2)e^{-2\pi m_2n_2/(32t)},
\end{align*}
where $t>0$ and the modular involution~\eqref{k01} was used.

Now,
\begin{align*}
L(E_{32},2)
&=L(f_{32},2)
=\int_0^1f_{32}\log q\,\frac{\d q}q
=-4\pi^2\int_0^\infty f_{32}(it)t\,\d t
\displaybreak[2]\\
&=-2\pi^2\int_0^\infty\sum_{m_1,n_1,m_2,n_2\ge1}a(m_1)b(n_1)b(m_2)a(n_2)
\\ &\qquad\qquad\times
\exp\biggl(-2\pi\biggl(m_1n_1t+\frac{m_2n_2}{32t}\biggr)\biggr)\d t
\displaybreak[2]\\
&=-2\pi^2\sum_{m_1,n_1,m_2,n_2\ge1}a(m_1)b(n_1)b(m_2)a(n_2)
\\ &\qquad\times
\int_0^\infty\exp\biggl(-2\pi\biggl(m_1n_1t+\frac{m_2n_2}{32t}\biggr)\biggr)\d t.
\end{align*}
Here comes the crucial transformation of purely analytical origin: we make
the change of variable $t=n_2u/n_1$. It does not change the form of the integrand
but affects the differential, and we obtain
\begin{align*}
L(E_{32},2)
&=-2\pi^2\sum_{m_1,n_1,m_2,n_2\ge1}\frac{a(m_1)b(n_1)b(m_2)a(n_2)n_2}{n_1}
\\ &\qquad\times
\int_0^\infty\exp\biggl(-2\pi\biggl(m_1n_2u+\frac{m_2n_1}{32u}\biggr)\biggr)\d u
\displaybreak[2]\\
&=-2\pi^2\int_0^\infty
\sum_{m_1,n_2\ge1} a(m_1)a(n_2)n_2 e^{-2\pi m_1n_2u}
\\ &\qquad\qquad\times
\sum_{m_2,n_1\ge1} \frac{b(m_2)b(n_1)}{n_1} e^{-2\pi m_2n_1/(32u)}\d u.
\end{align*}
What are the resulting series in the product? The first one corresponds to
\begin{equation*}
\sum_{m,n\ge1}a(m)a(n)n\,q^{mn}
=\sum_{m,n\ge1}\biggl(\frac{-4}{mn}\biggr)n\,q^{mn}
=\sum_{n\ge1}n\biggl(\frac{-4}n\biggr)\frac{nq^n}{1+q^{2n}}
=\frac{\eta_2^4\eta_8^4}{\eta_4^4},
\end{equation*}
while the second one is
\begin{align*}
\sum_{m,n\ge1}\frac{b(m)b(n)}n\,q^{mn}
&=\sum_{m,n\ge1}\frac{q^{mn}}n-\frac{q^{(2m)n}}n-\frac{q^{m(2n)}}{2n}+\frac{q^{(2m)(2n)}}{2n}
\\
&=\frac12\sum_{m,n\ge1}\frac{2q^{mn}-3q^{2mn}+q^{4mn}}n
\\
&=-\frac12\,\log\prod_{m\ge1}\frac{(1-q^m)^2(1-q^{4m})}{(1-q^{2m})^3}
=-\frac12\,\log\frac{\eta_1^2\eta_4}{\eta_2^3},
\end{align*}
hence
\begin{equation*}
L(E_{32},2)
=\pi^2\int_0^\infty\frac{\eta_2^4\eta_8^4}{\eta_4^4}\bigg|_{\tau=iu}
\cdot\log\frac{\eta_1^2\eta_4}{\eta_2^3}\bigg|_{\tau=i/(32u)}\d u.
\end{equation*}
Applying the involution \eqref{k01} to the eta quotient under the logarithm sign
we obtain
\begin{equation*}
L(E_{32},2)
=\pi^2\int_0^\infty\frac{\eta_2^4\eta_8^4}{\eta_4^4}
\,\log\frac{\sqrt2\eta_8\eta_{32}^2}{\eta_{16}^3}\bigg|_{\tau=iu}\d u.
\end{equation*}

Now comes the modular magic: assisted with Ramanujan's knowledge~\cite{Be-Rama}
we choose a particular modular function $x(\tau):=\eta_2^4\eta_8^2/\eta_4^6$,
which ranges from 1 to~0 when $\tau\in(0,i\infty)$, and verify that
$$
\frac1{2\pi i}\,\frac{x\,\d x}{2\sqrt{1-x^4}}
=-\frac{\eta_2^4\eta_8^4}{\eta_4^4}\,\d\tau
\quad\text{and}\quad
\biggl(\frac{\sqrt2\eta_8\eta_{32}^2}{\eta_{16}^3}\biggr)^2
=\frac{1-x}{1+x}.
$$
Thus,
\begin{equation*}
L(E_{32},2)
=\frac\pi8\int_0^1\frac{x}{\sqrt{1-x^4}}\,\log\frac{1+x}{1-x}\,\d x.
\end{equation*}

The result is a \emph{period} in the sense of \cite{KZ01}, and as such
it can be compared with several other objects like values of
generalized hypergeometric functions or even Mahler measures~\cite{RV99,Rog11}.
This however involves a different set of routines which we do not touch here.

To summarize, in our evaluation of $L(E,2)=L(f,2)$
we first split $f(\tau)$ into a product of two Eisenstein series of weight~1
and at the end we arrive at a product of two Eisenstein(-like) series
$g_2(\tau)$ and $g_0(\tau)$ of weights 2 and~0, respectively, so that
$L(f,2)=cL(g_2g_0,1)$ for some algebraic constant~$c$. The latter object
is doomed to be a period as $g_0(\tau)$ is a logarithm of a modular function,
while $2\pi i\,g_2(\tau)\,\d\tau$ is, up to a modular function multiple,
the differential of a modular function, and finally any two modular
functions are tied up by an algebraic relation over~$\overline{\mathbb Q}$.

The method however can be formalized to even more general settings,
and it is this extension which we attempt to outline below.

\bigskip
For two \emph{bounded} sequences $a(m)$, $b(n)$,
we refer to an expression of the form
\begin{equation}
g_k(\tau)=a+\sum_{m,n\ge1}a(m)b(n)n^{k-1}q^{mn},
\qquad q:=e^{2\pi i\tau},
\label{k03}
\end{equation}
as to an Eisenstein-like series of weight $k$, especially in the case
when $g_k(\tau)$ is a modular form of certain level, that is, when it transforms
sufficiently `nice' under $\tau\mapsto-1/(N\tau)$ for some positive integer~$N$.
This automatically happens when $g_k(\tau)$ is indeed an Eisenstein series
(for example, when $a(m)=1$ and $b(n)$ is a Dirichlet character modulo~$N$
of designated parity, $b(-1)=(-1)^k$),
in which case $\wh g_k(\tau):=g_k(-1/(N\tau))(\sqrt{-N}\tau)^{-k}$ is again an Eisenstein series.
It is worth mentioning that the above notion has perfect sense in case $k\le0$ as well.
Indeed, modular units, or week modular forms of weight 0, that are
the logarithms of modular functions are examples of Eisenstein-like series $g_0(\tau)$.
Also, for $k\le0$ examples are given by Eichler integrals,
the $(1-k)$\,th $\tau$-derivatives of holomorphic Eisenstein series of weight $2-k$,
a consequence of the famous lemma of Hecke~\cite[Section~5]{W77}.

Suppose we are interested in the $L$-value $L(f,k_0)$ of a cusp form $f(\tau)$ of weight $k=k_1+k_2$
which can be represented as a product (in general, as a linear combination of several products)
of two Eisenstein(-like) series $g_{k_1}(\tau)$ and $\wh g_{k_2}(\tau)$,
where the first one vanishes at infinity ($a=g_{k_1}(i\infty)=0$ in~\eqref{k03})
and the second one vanishes at zero ($\wh g_{k_2}(i0)=0$). (The vanishing happens
because the product is a cusp form!) In reality, we need the series
$g_{k_2}(\tau):=\wh g_{k_2}(-1/(N\tau))(\sqrt{-N}\tau)^{-k_2}$ to be Eisenstein-like:
$$
g_{k_1}(\tau)=\sum_{m,n\ge1}a_1(m)b_1(n)n^{k_1-1}q^{mn}
\quad\text{and}\quad
g_{k_2}(\tau)=\sum_{m,n\ge1}a_2(m)b_2(n)n^{k_2-1}q^{mn}.
$$

We have
\begin{align*}
L(f,k_0)
&=L(g_{k_1}\wh g_{k_2},k_0)
=\frac1{(k_0-1)!}\int_0^1 g_{k_1}\wh g_{k_2}\log^{k_0-1}q\,\frac{\d q}q
\\
&=\frac{(-1)^{k_0-1}(2\pi)^{k_0}}{(k_0-1)!}
\int_0^\infty g_{k_1}(it)\wh g_{k_2}(it)t^{k_0-1}\,\d t
\displaybreak[2]\\
&=\frac{(-1)^{k_0-1}(2\pi)^{k_0}}{(k_0-1)!\,N^{k_2/2}}
\int_0^\infty g_{k_1}(it)g_{k_2}(i/(Nt))t^{k_0-k_2-1}\,\d t
\displaybreak[2]\\
&=\frac{(-1)^{k_0-1}(2\pi)^{k_0}}{(k_0-1)!\,N^{k_2/2}}
\int_{0}^\infty\sum_{m_1,n_1\ge1}a_1(m_1)b_1(n_1)n_1^{k_1-1}e^{-2\pi m_1n_1t}
\\ &\qquad\qquad\times
\sum_{m_2,n_2\ge1}a_2(m_2)b_2(n_2)n_2^{k_2-1}e^{-2\pi m_2n_2/(Nt)}t^{k_0-k_2-1}\d t
\displaybreak[2]\\
&=\frac{(-1)^{k_0-1}(2\pi)^{k_0}}{(k_0-1)!\,N^{k_2/2}}
\sum_{m_1,n_1,m_2,n_2\ge1}a_1(m_1)b_1(n_1)a_2(m_2)b_2(n_2)n_1^{k_1-1}n_2^{k_2-1}
\\ &\qquad\times
\int_{0}^\infty\exp\biggl(-2\pi\biggl(m_1n_1t+\frac{m_2n_2}{Nt}\biggr)\biggr)t^{k_0-k_2-1}\d t;
\end{align*}
the interchange of integration and summation is legitimate because of the
exponential decrease of the integrand at the endpoints.
After performing the change of variable $t=n_2u/n_1$
and interchanging back summation and integration we obtain
\begin{align*}
L(f,k_0)
&=\frac{(-1)^{k_0-1}(2\pi)^{k_0}}{(k_0-1)!\,N^{k_2/2}}
\sum_{m_1,n_1,m_2,n_2\ge1}a_1(m_1)b_1(n_1)a_2(m_2)b_2(n_2)n_1^{k_1+k_2-k_0-1}n_2^{k_0-1}
\\ &\qquad\times
\int_{0}^\infty\exp\biggl(-2\pi\biggl(m_1n_2u+\frac{m_2n_1}{Nu}\biggr)\biggr)u^{k_0-k_2-1}\d u
\displaybreak[2]\\
&=\frac{(-1)^{k_0-1}(2\pi)^{k_0}}{(k_0-1)!\,N^{k_2/2}}
\int_{0}^\infty\sum_{m_1,n_2\ge1}a_1(m_1)b_2(n_2)n_2^{k_0-1}e^{-2\pi m_1n_2u}
\\ &\qquad\qquad\times
\sum_{m_2,n_1\ge1}a_2(m_2)b_1(n_1)n_1^{k_1+k_2-k_0-1}e^{-2\pi m_2n_1/(Nu)}u^{k_0-k_2-1}\d u
\displaybreak[2]\\
&=\frac{(-1)^{k_0-1}(2\pi)^{k_0}}{(k_0-1)!\,N^{k_2/2}}
\int_0^\infty g_{k_0}(iu)g_{k_1+k_2-k_0}(i/(Nu))u^{k_0-k_2-1}\,\d u.
\end{align*}
Assuming a modular transformation of the Eisenstein-like series
$g_{k_1+k_2-k_0}(\tau)$ under $\tau\mapsto-1/(N\tau)$, we can realize the
resulting integral as $c\pi^{k_0-k_1}L(g_{k_0}\wh g_{k_1+k_2-k_0},k_1)$,
where $c$~is algebraic
(plus extra terms when $g_{k_1+k_2-k_0}(\tau)$ is an Eichler integral).
Alternatively, if $g_{k_0}(\tau)$ transforms under the involution,
we perform the transformation and switch to the variable $v=1/(Nu)$ to
arrive at $c\pi^{k_0-k_1}L(\wh g_{k_0}g_{k_1+k_2-k_0},k_1)$. In both cases
we obtain an identity which relates the starting $L$-value $L(f,k_0)$
to a different `$L$-value' of a modular-like object of the same weight.

The case $k_1=k_2=1$ and $k_0=2$, discussed in \cite{RZ10,RZ11} and in our
example above, allows one to reduce the $L$-values to periods.
We wonder whether there are some other naturally `interesting'
(perhaps, known) examples. Already the case $k_0=k_1=k_2=1$ may be of interest,
when the $L$-values for two different cusp forms of weight~2 can be potentially compared.

\bigskip
\noindent
\textbf{Acknowledgements.}
I am thankful to the organizers of the RIMS conference ``Analytic number
theory\,---\,related multiple aspects of arithmetic functions'' (Kyoto University, Japan,
October~31--November~2, 2011) represented by Takumi Noda for invitation to
give a talk at the meeting. Special thanks go to my host Yasuo Ohno and
his team from the Kinki University (Osaka); they made my stay in Japan
both culturally and scientifically enjoyable.

I am indebted to Anton Mellit and Mat Rogers for fruitful conversations
on the subject, and to Don Zagier for his encouragement to isolate the transformation
part from \cite{RZ10,RZ11}.

\end{document}